\documentclass[12pt]{article}
\usepackage{cite,amsfonts,amssymb}
\usepackage{latexsym}
\newtheorem{prop}{Proposition}

\newtheorem{theor}{Theorem}

\newtheorem{cor}{Corollary}

\newtheorem{defin}{Definition}

\newcommand{\bgproof}{\noindent {\bf Proof} \hspace{2mm}}
\newcommand{\edproof}{\hfill $\blacksquare$ \vspace{3mm}}

\title{Temporal and spatial regularity of solutions to stochastic Volterra equations of convolution type}

\author{{\large\sf Anna Karczewska }\\[2mm]
  \normalsize\it
 Faculty of Mathematics, Computer Science and Econometrics\\ \normalsize\it
 University of Zielona G\'ora\\ \normalsize\it
 ul. Szafrana 4a, 65-246 Zielona G\'ora, Poland\\ \normalsize\it
 e-mail: A.Karczewska@wmie.uz.zgora.pl\\[2mm]
 }

\begin{document}

\maketitle

\def\thefootnote{}

\begin{abstract}
In the paper regularity of solutions to stochastic Volterra equations in a separable Hilbert space is studied.
Sufficient conditions for the temporal and spatial regularity of stochastic convolutions corresponding to the equations 
under consideration are provided.  The results obtained generalize some well-known regularity results for solutions 
to stochastic differential equations. The paper is a continuation of  previous author's papers concerning stochastic 
Volterra equations.
\end{abstract}

\section{Introduction}\label{sSW1}

Assume that 
$(\Omega,\mathcal{F},(\mathcal{F}_t)_{t\geq 0},P)$ is a stochastic basis.
We study stochastic linear Volterra equations in  a separable Hilbert space
 $(H, |\cdot|_H  )$ of the form
\begin{equation} \label{eSW1}
X(t) = X_0 + \int_0^t a(t-\tau)\,AX(\tau)d\tau + W(t)\;, \quad\quad t\geq 0\;.
\end{equation}
In (\ref{sSW1}),
 $X_0$ is an $ H$-valued, $\mathcal{F}_0$-measurable random variable, a scalar kernel function $a\in L^1_{\mathrm{loc}}(\mathbb{R}_+;\mathbb{R})$ 
and $A$ is a closed unbounded linear operator in $H$ with a dense domain
$D(A)$. In the domain $D(A)$ we introduce the graph norm $|\cdot
|_{D(A)}$ of $A$, i.e.\ $|h|_{D(A)}:=(|h|_H^2+|Ah|_H^2)^{1/2}$.
In the paper we study the equation (\ref{eSW1}) which is driven by a 
$Q$-Wiener process $W$ with $\mathrm{Tr}\,Q <+\infty$. 

Let us emphasize that the equation (\ref{eSW1}) includes a big class of equations, for instance heat and wave equations and integrodifferential equations. Moreover, it is an abstract version of several technical problems.

Regularity of solutions to equations is very important and plays a prominent role in the study of evolution equations, particularly stochastic ones, see e.g. \linebreak \cite{DaPrZa, PeZa}.  In the paper 
we expose sufficient conditions for time and space regularity of stochastic convolution, which is a solution 
to the stochastic equation~(\ref{eSW1}). In our study we use the so-called resolvent approach to 
the equation~(\ref{eSW1}), which is a natural extension of the semigroup approach used to stochastic 
differential equations. We note that the results obtained in the paper are an extension of the regularity results 
for solutions to stochastic differential equations given, e.g., in \cite{DaPrZa}.

Let us note that our case considered in this paper is more difficult than the semigroup case connected with stochastic differential equations. The convolution appearing in~(\ref{eSW1}) causes several serious problems and difficulties; some of them have been discussed in~\cite{Ka10}. The most big problem comes from the fact that the solution operator (or the resolvent operator) corresponding to~(\ref{eSW1}) and recalled below, does not create, in general, any semigroup. In consequence, powerful semigroup tools can not be used to the study of the equation~(\ref{eSW1}). Particularly, the factorization method due to G.~Da Prato, S.~Kwapie\'n and J.~Zabczyk, provided in \cite{DPKZ87}, is not available in the resolvent case. However, some regularity results for solutions to~(\ref{eSW1}) are known, e.g., those due to Ph.~Clement and G.~Da Prato \cite{ClDP96}.

Our approach to regularity of solution to~(\ref{eSW1}) is different. We rewrite the stochastic convolution associated with~(\ref{eSW1}), which is a mild and/or weak solution to~(\ref{eSW1}). Next, we use the new form (see formula (\ref{eR1}) below) of the stochastic convolution and join it with appropriate Cauchy problem (see formulas (\ref{eco10}) and (\ref{eco10a}) below). Then we use regularity results of some Cauchy problems considered in monographs \cite{Lun95,Lun09,Pa}, adapting them for our purposes.

The paper is organized as follows. Section 2 contains an appropriate 
deterministic background while section 3 provides the main definitions and concepts 
used in the paper. In section 4 we give sufficient conditions for stochastic convolution
to be regular solution to the equation (\ref{eSW1}). The final section 5 consists of proofs of the main results.

\section{Deterministic Volterra equation}\label{sSW2}

The following equation 

\begin{equation} \label{eSW1d}
u(t) = \int_0^t a(t-\tau)\,Au(\tau)d\tau +f(t), \quad t\geq 0,
\end{equation} 
is a deterministic counterpart of  (\ref{eSW1}) in the space $H$.
In the equation  (\ref{eSW1d}), the operator $A$ and the
kernel function $a$ are the same as previously and $f$ is a continuous $H$--valued function.

In the paper we assume that the equation  (\ref{eSW1d}) is well-posed, that is,  the family  $S(t),~t\geq 0$, 
of the resolvent operators corresponding to the Volterra equation (\ref{eSW1d}) exists and is defined as follows.

\begin{defin}\label{dSW1} (See, e.g., \cite{Pr2}.)~ 
A family $(S(t))_{t\geq 0}$ of bounded linear operators in $H$ is
called {\tt resolvent} for the equation  (\ref{eSW1d}) if the following
conditions are satisfied:
\begin{enumerate}
\item $S(t)$ is strongly continuous on $\mathbb{R}_+$ and $S(0)=I$;
\item $S(t)$ commutes with the operator $A$, that is, $S(t)(D(A))\subset
D(A)$ and $AS(t)x=S(t)Ax$ for all $x\in D(A)$ and $t\geq 0$;
\item the following {\tt resolvent equation} 
\begin{equation} \label{eSW2}
S(t)x = x + \int_0^t a(t-\tau) AS(\tau)x d\tau
\end{equation}
holds for all $x\in D(A),~t\geq 0$.
\end{enumerate}
\end{defin}

Because the class of locally integrable scalar kernel functions is too wide for obtaining good enough results concerning 
regularity of solutions to the equation (\ref{eSW1}), we shall consider the case when the kernel function 
is completely positive understood in the following sense.

\begin{defin} \label{dSW2}
We say that function $a\in L^1([0,T];\mathbb{R})$ is {\tt completely positive}
on $[0,T]$, $T<+\infty$, if for any $\mu\geq 0$, the solutions of the equations
\begin{equation}
s(t)\! +\! \mu \!\int_0^t\!  a(t-\tau)s(\tau)\,d\tau =1 \quad \mbox{and} \quad r(t) \!+\! \mu\! \int_0^t\!  a(t-\tau)r(\tau)\,d\tau = a(t)
\end{equation}
are nonnegative on $[0,T]$.
\end{defin}

There are several examples of completely positive kernels, e.g.:
\begin{enumerate}
\item $a(t)= t^{\alpha-1}/\Gamma(\alpha), \quad \alpha\in (0,1],\quad t>0$;
\item $a(t)= e^{-t}, \quad t\ge 0 $.
\end{enumerate}

\noindent The class of
completely positive kernels appears in the theory of
viscoelasticity. Several properties and examples of such kernels
can be found in \cite[Section 4.2]{Pr2}.

In the paper we shall use the following theorem providing convergence of the resolvents for the equation (\ref{eSW1d}).

\begin{theor} \label{pSW2} (\cite[Theorem 4]{KaLiJMAA})~
Let $A$ be the generator of $C_0$-semigroup  in $H$ and suppose
the kernel function $a$ is completely positive. Then the pair $(A,a)$
admits an exponentially bounded resolvent $S(t), ~t\ge 0$. Moreover, there
exist bounded operators $A_n$ such that $(A_n,a)$ admit resolvent
families $S_n(t), ~t\ge 0$, satisfying $ ||S_n(t) || \leq Me^{w_0 t}~ (M\geq
1,~w_0\geq 0)$ for all $t\geq 0$ and
\begin{equation} \label{eSW4}
S_n(t)x \to S(t)x \quad \mbox{as} \quad n\to +\infty,
\end{equation}
for all $x \in H,\; t\geq 0.$ Additionally, the convergence is
uniform in $t$ on every compact subset of $ \mathbb{R}_+$.
\end{theor}

In the above theorem and in the sequel, we use the Yosida approximation of the operator $A$ and we denote it by
\begin{equation} \label{eSW5}
A_n := n AR(n,A) = n^2 R(n,A) - nI, \qquad n\in \mathbb{N}.
\end{equation}

The analogous results like Theorem \ref{pSW2} hold in other cases; for more details concerning the approximation results we refer to \cite{Ka07,KaLiJEE,KaLiJMAA}.

\section{Probabilistic background} \label{sSW3}

Let  $Q\in L(H)$ be a linear bounded symmetric nonnegative operator acting in the space $H$. We recall that the Wiener process $W$  with the covariance operator $Q$  is assumed to be a genuine one, that is,  $\mathrm{Tr} Q<+\infty$. 

In the paper we shall consider the following types of solutions  to the equation~(\ref{eSW1}).

Let $A^*$ denote the adjoint of $A$ with a dense domain
$D(A^*)\subset H$ equipped with the graph norm $|\cdot |_{D(A^*)}$.

\begin{defin} \label{dSW5}
 An $H$-valued predictable process $X(t),~t\in
[0,T]$, is said to be a {\tt weak solution} to (\ref{eSW1}), if
$P(\int_0^t|a(t-\tau)X(\tau)|_H d\tau<+\infty)=1$ and if for all
$\xi\in D(A^*)$ and all $t\in [0,T]$ the following equation holds
$$
\langle X(t),\xi\rangle_H = \langle X_0,\xi\rangle_H + \langle
\int_0^t a(t-\tau)X(\tau)\,d\tau, A^*\xi\rangle_H + \langle
W(t),\xi\rangle_H, ~P\mathrm{-a.s.}
$$
\end{defin}

As we have already written we assume that (\ref{eSW1d}) admits a resolvent
family $S(t)$, $t \geq 0.$  Then we are able to use the resolvent approach to 
both equations (\ref{eSW1d}) and  (\ref{eSW1}). So, we can introduce the following definition.

\begin{defin} \label{dSW6}
An $H$-valued predictable process $X(t),~t\in
[0,T]$, is said to be a {\tt mild solution} to the stochastic
Volterra equation (\ref{eSW1}), if
 for arbitrary $t\in [0,T]$,
\begin{equation}\label{eSW9}
X(t) = S(t)X_0 + \int_0^t S(t-\tau)dW(\tau), \quad
P-a.s.
\end{equation}
where $S(t), ~t\ge 0$, is the resolvent for the equation (\ref{eSW1d}).
\end{defin}

 We introduce the stochastic convolution, which is the main part of the mild solution 
\begin{equation} \label{eSW18a}
W^S(t) := \int_0^t S(t-\tau)dW(\tau), \qquad t\ge 0.
\end{equation}

Let us recall some  results concerning the convolution
$W^S (t), ~t\ge 0$, proved in \cite{Ka}.

\begin{prop}\label{pr2proper} (See \cite[Propositions~3 and 4]{Ka}.) 
Assume that (\ref{eSW1d}) admits the resolvent operators $S(t),~t\geq
0$. Then 
the process $W^S(t),~t\geq 0$, given by (\ref{eSW18a}) has a
predictable version. Moreover, the process
$W^S(t),~t\geq 0$,  has square
integrable trajectories.
\end{prop}

In some cases weak solution of the equation (\ref{eSW1}) coincides
with mild solution of (\ref{eSW1}); for more information we can refer to \cite{Ka,Ka07}.

\begin{prop} \label{pr4proper} (See \cite[Proposition~5]{Ka}.)
Let $a\in BV(\mathbb{R}_+;\mathbb{R})$ and suppose that (\ref{eSW1d}) admits
a resolvent family $S \in C^1(0,\infty; L(H)).$ Let $X$ be a
predictable process with integrable trajectories. Assume that $X$
has a version such that $P(X(t)\in D(A))=1$ for almost all $t\in
[0,T]$. If for any $t\in [0,T]$ and $\xi\in
D(A^*)$,
\begin{eqnarray}\label{deq9}
 \langle X(t),\xi\rangle_H = \langle X_0,\xi\rangle_H  &+& 
 \int_0^t \langle
 a(t-\tau)X(\tau),A^*\xi\rangle_H d\tau \\ &+&
 \langle W(t),\xi\rangle_H, ~P-a.s., \nonumber
\end{eqnarray}
then
\begin{equation}\label{deq9a}
X(t) = S(t)X_0 +
 \int_0^t S(t-\tau) dW(\tau), \quad t\in[0,T].
\end{equation}
\end{prop}

\begin{prop} \label{pr5prop} (See \cite[Proposition~6]{Ka}.)
Assume that $A$ is a closed linear unbounded operator with the
dense domain $D(A)$, $a\in L_\mathrm{loc}^1(\mathbb{R}_+;\mathbb{R})$ and
$S(t)$, \linebreak $t\ge 0$, are resolvent operators for the equation
(\ref{eSW1d}). Then the
stochastic convolution $W^S$ fulfills the equation (\ref{deq9})
with $X_0\equiv 0$.
\end{prop}

\begin{cor} \label{c1cor} (Compare \cite[Corollary 3.10]{Ka07}.)
By Proposition \ref{pr5prop}, if $A$ is a linear bounded operator, we have
\begin{equation} \label{deq16}
 W^S(t) = \int_0^t a(t-\tau)AW^S(\tau)d\tau
 + W(t), \quad t\in [0,T].
\end{equation}
\end{cor}

Let us note that the above formula (\ref{deq16}) indicates continuity of trajectories of the convolution $W^S(t)$, $ t\in [0,T]$, in the case when $A$ is a bounded operator. In our paper we shall show continuity of trajectories in a more general case when $A$ is an unbounded operator. Additionally, we will prove some regularity results for the convolution (\ref{eSW18a}).

\section{Time and space regularity}\label{SS4} 

In this section we study the temporal and spatial regularity of the stochastic convolution $W^S$, 
defined by (\ref{eSW18a}), that is
$$W^S(t) = \int_0^t S(t-\tau)dW(\tau)\,, \qquad t\in [0,T], \qquad  T<+\infty .$$

\begin{theor} \label{R1}
Assume that the operator $A$ in the equation (\ref{eSW1}) is the generator of a $C_0$-semigroup $T(t), t\in [0,T]$. Let the kernel function $a(t), t\in [0,T]$, 
be completely positive and, additionally, 
$\dot{a}\in L^1_\mathrm{loc} ([0,T];\mathbb{R})$. 
Let $S(t)$, $W(t)$ and $W^S(t), t\in [0,T],$
 be like above. Then the following formula holds
\begin{equation} \label{eR1} \hspace{-4ex}
W^S(t) =c A\int_0^t T(t-\tau) \left[ \int_0^\tau \dot{a}(\tau-\sigma) W^S(\sigma) d\sigma+c W(\tau) \right]d\tau +W(t), 
\end{equation} $t\in [0,T]$,
where $c=a(0)$ is a constant.
\end{theor}

Examples of functions, fulfilling assumptions of Theorem \ref{R1}, are:
\begin{enumerate}
\item $a(t)= e^{-t}, \quad t\in [0,T] $;
\item $a(t)= t^{\alpha-1}/\Gamma(\alpha), \quad \alpha\in (0,1],\quad t\in [\varepsilon,T],\quad \varepsilon > 0$.
\end{enumerate}

\noindent
{\bf Comment~} 
One can see that the above formula (\ref{eR1}) is more complicated than the analogous one for the 
semigroup case, see, {\em e.g.} \cite[Chapter 5]{DaPrZa}. In fact, if $a(t)=1$, the formula (\ref{eR1}) reduces to 
$$ W^S(t) = A\int_0^t T(t-s) W(s)ds + W(t), \qquad t\in [0,T]\,. $$
\\

For our convenience we will assume in the sequel that $c\equiv a(0)=1$.

\begin{cor} \label{co10}
 Suppose that the assumptions of Theorem \ref{R1} hold and $a(0)\!=\!1$. Define the process 
\begin{equation} \label{eco10}
 Y(t) := \int_0^t T(t-s)\left[ \widetilde{W}(s)+W(s)\right] ds, \qquad t\in [0,T]\,,
\end{equation}
where 
$$  \widetilde{ W}(s):= \int_0^s \dot{a}(s-\sigma) \,W^S(\sigma) \,d\sigma, \qquad s\in[0,T]\,. $$
Then
\begin{equation} \label{eco10x}
 W^S(t) = A\,Y(t) + W(t), \qquad t\in [0,T]. 
\end{equation}

 Additionally, $Y$ belongs to $C^1([0,T];D(A)),~ \mathbb{P}-a.s.$, and
\begin{equation} \label{eco10a}
 \frac{dY(t)}{dt} = A\,Y(t) + \left[ \widetilde{W}(t)+W(t)\right], \qquad t\in[0,T]\,.
\end{equation}
\end{cor}

From now we assume that the space $H$ is a complex separable Hilbert space and the operator $A$ is the generator of a strongly continuous  analytic semigroup.

Because both processes $Y$ and $W$ have continuous trajectories, directly from the formula (\ref{eco10x}) we can deduce the following result.

\begin{theor} \label{R2}
Let us take the same assumptions like in Theorem \ref{R1} and additionally let the semigroup be analytic. Then the stochastic convolution $W^S(t)$, $t\ge 0$,
has continuous trajectories.
\end{theor}

\noindent
{\bf Comment~} Let us note that Theorem  \ref{R1}, Corollary  \ref{co10} and Theorem \ref{R2} can be formulated for other cases when the convergence  (\ref{eSW4}) of the resolvents for the equation  (\ref{eSW1d}) holds. For instance in the cases described in the  papers \cite{Ka07,KaLiPAMS,KaLiJEE}. \\

Unfortunately, we are not able to obtain H\"olderianity of trajectories of the stochastic convolution~(\ref{eSW18a}). In contrary to the semigroup case, the formula (\ref{eR1}) for the stochastic convolution is more complicated and contains the term $\widetilde{W}$ (see formula (\ref{eco10})). In consequence, the nonhomogeneous Cauchy problem (\ref{eco10a}) contains two functions $W$ and $\widetilde{W}$. The first one has  H\"older-continuous trajectories, but the second one has not such trajectories.

There are available some temporal regularity results for stochastic Volterra equations obtained under different assumptions, see e.g., \cite{ClDP96} and \cite{KaLiPAMS}.

In order to study next properties of the stochastic convolution (\ref{eSW18a}) we shall introduce some scales of subspaces of $H$. We will use some results for analytic semigroups generated by sectorial operators and the Cauchy problems associated with such semigroups.

\begin{defin} \label{sec_op} (See, e.g., \cite{Lun95,Lun09}.)~
A linear operator $A : D(A)\subset H \to H$ is called {\tt sectorial} if there are constants $\omega\in\mathbb{R}$, $\theta\in (\pi/2,\pi)$, $M>0$ such that
\begin{equation}  \label{esec}
\left\{ \begin{array}{ll} ~(1)~ & \rho(A) \supset S_{\theta,\omega} := \{\lambda\in \mathbb{C}: \lambda\ne \omega,~ |\arg(\lambda - \omega)| <\theta \}\,, \\ 
~(2)~ & ||R(\lambda,A) ||_{L(H)} \le \displaystyle\frac{M}{ |\lambda - \omega|} \quad \forall \, \lambda\in S_{\theta,\omega} \,.  \end{array} \right.
\end{equation}
\end{defin}

Assume that the operator $A$ is sectorial, then we can consider an analytic semigroup $T(t),~ t\ge 0$, defined in $H$ through the Dunford integral. For details on sectorial operators and analytic semigroups
we refer to the great monographs \cite{Bensou,EnNa} and \cite{Lun95,Lun09,Pa}.
If the semigroup $T(t),~ t\ge 0$, is assumed to be of negative type, we can define the  fractional powers of operators as follows. For any $\gamma\in(0,1)$, we set
$$ (-A)^{-\gamma} x \stackrel{\rm{df}}{=} \frac{1}{2\pi i} \int_{\Gamma_{r,\theta}} (-\lambda)^{-\gamma} R(\lambda,A)\,x\,d\lambda\,, \qquad x\in H, $$
where $\lambda$ is a complex number and $R(\lambda,A):=(\lambda I-A)^{-1}$ is the resolvent operator of $A$. The curve $\Gamma_{r,\theta}$ is defined as follows. For $\theta\in (\pi/2,\pi)$,\linebreak $r>0$, $\Gamma_{r,\theta}= -\Gamma_{r,\theta}^{1}-\Gamma_{r,\theta}^{2}+\Gamma_{r,\theta}^{3}$, where $\Gamma_{r,\theta}^{1},\Gamma_{r,\theta}^{3}$ are the half lines parametrized respectively by $z = \xi e^{i\theta},~ z = \xi e^{-i\theta},~ \xi\ge r$, and $\Gamma_{r,\theta}^{2}$ is the arc of circle prametrized by $z = r e^{i\eta}, ~ -\theta\le\eta\le\theta$.

By $(-A)^{\gamma}$ we shall denote the inverse of the operator $(-A)^{-\gamma}$ and by $D((-A)^{\gamma})$ its domain. In the literature the operators $(-A)^{\gamma}$ are called the fractional powers of the operator $-A$. They have very interesting and useful properties, see e.g.\ \cite{Lun95,Lun09,Pa}. Among others, there exist contants $M_k, M_{k\gamma}$, $k=0,1$, $\gamma\in(0,1)$, such that
\begin{equation} \label{nn1}
|| A^k\,T(t)|| \le M_k\,t^{-k}, \qquad t\ge 0,
\end{equation}
and 
\begin{equation} \label{nn2}
|| (-A)^\gamma A^k\,T(t)|| \le M_{k\gamma}\, t^{-k-\gamma}, \qquad t\ge 0.
\end{equation}

In the remaining part of the paper we shall use the following properties of powers of operators.

\begin{theor} \label{Lun4.6} ($\!$\cite[Theorem 4.6]{Lun09})~
Let $\alpha,\beta\in\mathbb{C}$ such that {\rm $\mbox{Re}\;\beta<\mbox{Re}\;\alpha $}. Then $D(A^\alpha)\subset D(A^\beta)$, and for every $x\in D(A^\alpha)$
$~A^\beta x = A^{\beta-\alpha} A^\alpha x.$~
Moreover, for each $x\in D(A^\alpha)$, $A^\beta x\in D(A^{\beta-\alpha})$ and
$~ A^{\alpha-\beta} A^\beta x = A^\alpha x.$
Conversly, if $x\in  D(A^{\beta})$ and $A^\beta x\in D(A^{\alpha-\beta})$, then $x\in D(A^\alpha)$ and again $ A^{\alpha-\beta} A^\beta x = A^\alpha x$.
\end{theor}

Now, we are able to formulate spatial regularity for the convolution  (\ref{eSW18a}).
\begin{theor} \label{R3}
For all $t\in [0,T]$ and $\gamma\in(0,1)$, the convolution $W^S(t)$ belongs to the space $D((-A)^\gamma)$, a.s., provided the process $W$ belongs to the domain $D((-A)^\gamma)$. 
Additionally, $(-A)^\gamma W^S$ is a Gaussian process and
 for all \linebreak $\gamma\in(0,1)$ trajectories of the process $(-A)^\gamma W^S$ are continuous.
\end{theor}

Following A.~Lunardi \cite{Lun95,Lun09} we introduce the interpolation spaces.

\begin{defin}  \label{dLu1}
~For ~$\gamma\in (0,1)$, $1\le p\le +\infty$, we denote $D_A(\gamma,p):=(H,D(A))_{\gamma,p}$, where the pair $(H,D(A))_{\gamma,p}$ is the interpolation spaces according to e.g.\ Definition 1.2 in \cite{Lun09}.
\end{defin}

\begin{defin}  \label{dLu2}
~For ~$\gamma\in (0,1)$, $1\le p\le +\infty$, we define $D_A(\gamma+k,p):= \{ x\in D(A^k):A^k x\in D(A)_{\gamma,p}\}$, with the norm $|x|_{D_A(\gamma+k,p)}:=|x|_H+|A^k x|_{D_A(\gamma,p)}$. 
\end{defin}

Let us note that $D_A(\gamma+k,p)$ is the domain of the part of $A^k$ in $D_A(\gamma,p)$.

Now, we are able to formulate further regularity results for the stochastic convolution $W^S$ defined by (\ref{eSW18a}) in terms of interpolation spaces.

The formula (\ref{eco10a}) enables us to study regularity of the process $Y$ using properties of deterministic Cauchy problem.

By \cite[Proposition 6.6]{Lun09} and from square integrability on $[0,T]$ of trajectories of the processes $W$ and $\widetilde{W}$ we have the following result.

\begin{theor} \label{R5}
~For every $\gamma\in (0,1)$, trajectories of the process $Y$  belong to the space $L^2((0,T);D_A(\gamma,2))\cap W^{\gamma,2}((0,T);H)$.  Consequently, it belongs to the space $ W^{\gamma-\theta,2}((0,T);D_A(\theta,2))$ for every $\theta<\gamma<1$ and to $C([0,T];D_A(\gamma,2))$ for every $\gamma\in(0,\frac{1}{2})$. Moreover, there is a constant $M$ independent of the processes $W$ and $\widetilde{W}$, such that 
$$ |Y|_{W^{\gamma,2}((0,T);H)}+|Y|_{L^2((0,T):D_A(\gamma,2))} \le M\, |W+\widetilde{W}|_{L^2((0,T);H)}. $$
\end{theor}

The next result is formulated without using interpolation of powers of operators. We shall use so-called Rademacher bounded operator.
\begin{defin}  \label{dRad}
~A subset of $L(H)$ is $\mathbb{R}$-bounded (Rademacher bounded) if there is $C>0$ such that for all $n\in\mathbb{N}$, $L_1,\ldots ,L_{n}\in L$, $x_1,\ldots ,x_{n}\in H$ we have
$$ \left|\left| \sum_{k=1}^n r_k T_k x_k\right|\right|_{L^2((0,1:H))} \le \left|\left| \sum_{k=1}^n r_k x_k\right|\right|_{L^2((0,1:H))},$$
where $r_k(t):= sign(\sin 2^k\pi t)$ are Rademacher functions in $(0,1)$.
\end{defin}

Using \cite[Theorem 6.11]{Lun09} we can formulate the following theorem.

\begin{theor} \label{R6}
~Assume that the operator $A: D(A)\subset H\to H$ is sectorial and the semigroup generated by $A$ is strongly continuous. Then the following two statements are equivalent: 
\begin{enumerate}
\item trajectories of the process $Y$ belong to the space \linebreak
$\displaystyle W^{1,2} ((0,T);H)\cap L^2((0,T);D(A))$;
\item there are $\omega\in \mathbb{R}$, $\theta\ge \frac{\pi}{2}$, such that the family of operators \linebreak $\{(\lambda-\omega)R(\lambda,A):\lambda\in S_{\theta,\omega}\}$ is $\mathcal{R}$-bounded.
\end{enumerate}
\end{theor}

Let us note, see e.g.\ remark on page 171 in \cite{Lun09}, that in Hilbert space every sectorial operator satisfies the condition (2) in Theorem \ref{R6}. So, we can formulate the following conclusion.

\begin{cor} \label{R7}
~Assume that $A: D(A)\subset H\to H$ is a sectorial operator and the semigroup generated by $A$ is strongly continuous. Then trajectories of the process $Y$ belong to  $W^{1,2} ((0,T);H)\cap L^2((0,T);D(A))$.
\end{cor}

\section{Proofs of the main results}\label{SS5}

\bgproof  of Theorem~\ref{R1}.~\\
Because the formula (\ref{deq16}) holds for any bounded operator, then it holds for the Yosida approximation $A_n$ of the operator $A$, that is 
\begin{equation} \label{eR1a} 
W^{S_n}(t) =\int_0^t a(t-\tau)A_n W^{S_n}(\tau)\,d\tau +W(t), \quad t\in [0,T].
\end{equation} 
In the formula (\ref{eR1a}) , $S_n(t)$, $t\ge 0$, is the resolvent genereted by the pair $(A_n,a(t))$, $t\ge 0$, (see Theorem \ref{pSW2}) and 
\begin{equation} \label{eR1b} 
W^{S_n}(t) := \int_0^t S_n(t-\tau)\,dW(\tau)\,.
\end{equation}
Let us denote
\begin{equation} \label{eR1c}
Z_n(t) := \int_0^t a(t-\tau) W^{S_n}(\tau)\,d\tau, \quad t\in[0,T]\,.
\end{equation}
Then, from the Leibniz rule 
\begin{eqnarray} \label{eR1d}
(Z_n(t))' & = & \left(\int_0^t a(t-\tau) W^{S_n}(\tau)\,d\tau\right)' =\\
 & = & \int_0^t \dot{a}(t-\tau)\, W^{S_n}(\tau)\,d\tau + a(0)\,W^{S_n}(t)\,. \nonumber
\end{eqnarray}
From (\ref{eR1a}) and (\ref{eR1c}), we can write 
$$ W^{S_n}(t) = A_n \,Z_n(t) +W(t), \quad t\in[0,T]\,. $$
From (\ref{eR1d}), if $a(0)\ne 0$, we have 
$$ W^{S_n}(t) = \frac{1}{a(0)} \left[Z_n'(t) - \int_0^t \dot{a}(t-\tau)\, W^{S_n}(\tau)\,d\tau \right]\,.$$
So, we obtain
$$ Z_n'(t)  -  \int_0^t \dot{a}(t-\tau) W^{S_n}(\tau)\,d\tau =
 a(0)\, [A_n\,Z_n(t) +W(t)], \quad t\in [0,T]. $$
And next
$$ Z_n'(t) = a(0)\, A_n\,Z_n(t) + \int_0^t \dot{a}(t-\tau) W^{S_n}(\tau)\,d\tau +  a(0)\,W(t), \quad t\in [0,T]. $$ 
For simplicity, we set
$$ \widetilde{W}_n^S(t) := \int_0^t \dot{a}(t-\tau) W^{S_n}(\tau)\,d\tau, \quad t\in [0,T]. $$
So, we have 
$$ Z_n'(t) = c\, A_n\,Z_n(t) + \left[\widetilde{W}_n^S(t)+c\,W(t)\right], \quad 
\mbox{where} \quad c=a(0). $$
From the above formula 
$$ Z_n(t) = \int_0^t e^{c(t-\tau)A_n} \left[\widetilde{W}_n^S(\tau)+c\,W(\tau)\right]\,d\tau,\quad t\in [0,T]\,, $$
because $Z_n(0)=0$.

From the formula (\ref{eR1a}),
$$  W^{S_n}(t) =A_n\,Z_n(t) + W(t), \quad t\in [0,T]\,, \qquad \mbox{or} $$
$$  W^{S_n}(t) =  A J_n Z_n(t) + W(t), \qquad t\in [0,T],$$
where  $J_n := nR(n,A)$. 

So, we have 
$$  W^{S_n}(t) =AJ_n\int_0^t e^{c(t-\tau)A_n} \left[\widetilde{W}_n^S(\tau)+c\,W(\tau)\right]\,d\tau +W(t),\quad t\in [0,T]. $$
Basing on Theorem \ref{pSW2}, properties of Yosida approximation $A_n$ of the operator $A$, and the Lebesgue dominated convergence theorem, we have:
$$  \lim_{n\to\infty} J_n\,x = x,~ \quad \mbox{for~any} \quad x\in H;   \qquad$$  
$$ \lim_{n\to\infty} A_n\,x = A\,x, \quad \mbox{for~any} \quad x\in D(A); $$
$$ \lim_{n\to\infty} e^{tA_n}\,x = T(t)\,x, \quad \mbox{for~any} \quad x\in H; $$
$$\mbox{and} \quad \lim_{n\to\infty}  \sup_{t\in [0,T]} \mathbb{E}\left|W_n^S(t)-W^S(t)\right|_H^2= 0 . $$

Because the operator $A$ is closed, we can conclude that the integral \linebreak $\displaystyle \int_0^t T(t-\tau)\left[\widetilde{W}^S(\tau)+c\,W(\tau) \right]d\tau~$ belongs to the domain $D(A)$.

Hence, passing to the limit with $n\to +\infty$, we obtain 
$$  W^S(t) =c A\int_0^t T(t-\tau)\, \left[\widetilde{W}^S(\tau)+c\,W(\tau)\right]\,d\tau +W(t), \quad t\in [0,T], $$
where $$ \widetilde{W}^S(\tau) =\int_0^\tau \dot{a}(\tau-\sigma)W^S(\sigma)d\sigma\,,\quad \tau\in [0,T]\,. $$

In the case $a(0)=0$, from the formula (\ref{eR1a}), passing to the limit we would have only 
$$  W^S(t) = \int_0^t a(t-\tau)AW^S(\tau)d\tau + W(t), \quad t\in [0,T]\,, $$
that is, the formula like (\ref{deq16}). 
\edproof

\bgproof   of Corollary~\ref{co10}.~\\
From the definition (formula (\ref{eco10})) of the process $Y$ and properties of convolution, $Y\in C^1([0,T];D(A)),~ P-a.s.$  Next, from the Leibniz rule and property of semigroup we obtain
\begin{eqnarray*}
 \frac{dY(t)}{dt} &=& \int_0^t  \frac{dT(t-s)}{dt}\left[\widetilde{W}(s) +W(s) \right]ds + T(0) \left[\widetilde{W}(t) +W(t) \right] \\
 &=& A \int_0^t  T(t-s)\left[\widetilde{W}(s) +W(s) \right]ds + \left[\widetilde{W}(t) +W(t) \right] \\
&=&  A\,Y(t) + \left[\widetilde{W}(t) +W(t) \right], \qquad t\in[0,T].
\end{eqnarray*} 
\edproof

\bgproof   of Theorem~\ref{R3}.~\\
In our case we have no an appropriate representation formula for the resolvent like in the semigroup case, see, e.g., \cite[Chapter 2]{Pa}. In consequence, we have no analogous estimates for resolvents like (\ref{nn1})-(\ref{nn2}). In order to omit these difficulties, we have to use the formula (\ref{eco10x}) for the convolution $W^S$.

By (\ref{nn2}), the integral
$$ \int_0^t (-A)^\gamma AT(t-s)\left[\widetilde{W}(s) +W(s) \right]ds, \quad t\in  [0,T], \quad \gamma\in (0,1), $$ 
is well defined. Because $(-A)^\gamma$ is closed, it follows that $W^S\in D((-A)^\gamma)$.

Additionally, because $W^S\in D((-A)^\gamma)$, we have
$$ (-A)^\gamma W^S(t) =\! \int_0^t (-A)^\gamma AT(t-s)\left[\widetilde{W}(s) + W(s) \right]ds\; +\; (-A)^\gamma W(t),~~ t\in [0,T].
$$
So, the process $(-A)^\gamma W^S$ is well defined.

Gaussianity of the process comes from Gaussianity of the process $W^S$. Continuity of trajectories of the process $(-A)^\gamma W^S$ results from continuity of trajectories of the process $W^S$.
\edproof

\end{document}